\documentclass[reqno]{amsart}
\usepackage{amsmath,amsthm,amssymb,hyperref,graphicx,mathrsfs,mathtools}  

\newtheorem{theorem}{Theorem}
\newtheorem{corollary}{Corollary}
\newtheorem{lemma}{Lemma}

\newtheorem{remark}{Remark}

\allowdisplaybreaks

\begin{document}
\title{A note on generalized $L_p$ inequalities for the polar derivative of a  polynomial}
\author{N. A. Rather}
\author{DANISH RASHID BHAT}
\author{Tanveer Bhat}
\address{$^{1,2,3}$Department of Mathematics, University of Kashmir, Srinagar-190006, India}
\email{$^1$dr.narather@gmail.com, $^2$danishmath1904@gmail.com $^3$Tanveerbhat054@gmail.com}

\begin{abstract}
Let \( P(z) \) be a polynomial of degree \( n \) and $\alpha \in \mathbb{C}$. The polar derivative of \( P(z) \), denoted by \( D_\alpha P(z) \) and is defined by
$D_\alpha P(z): = nP(z) + (\alpha -z)P'(z)$. The polar derivative \( D_\alpha P(z) \) is a polynomial of degree at most \( n - 1 \) and it generalizes the ordinary derivative \( P'(z) \).
In this paper, we establish some \( L_p \) inequalities for the polar derivative of a polynomial with all its zeros located within a prescribed disk. Our results refine and generalize previously known findings.
\\
\smallskip
\newline
\noindent \textbf{Keywords:}  $L_p$ inequality, Polar Derivative, Gauss-Lucas Theorem.\\
\noindent \textbf{2010 Mathematics Subject Classification:} 26D10, 41A17, 30A10, 26D05
\end{abstract}

\maketitle

 \section{ Introduction}

Let $P(z)$ be a polynomial of degree $ n $. A well known result, known as Bernstein's inequality \cite{I}, provides an upper bound for the maximum modulus of $ P'(z)$ on the unit circle $|z| = 1$ in terms of the maximum modulus of $ P(z) $, which is given by\\
\begin{align*}
\displaystyle\max_{|z|=1} |P'(z)| \leq n \displaystyle\max_{|z|=1} |P(z)|.
\end{align*}
Equality holds if and only if all the zeros of $ P(z) $ are at the origin. For polynomials whose zeros lie in the unit disk $|z| \leq 1$, Turán \cite{TU} established a lower bound for the maximum modulus of $P'(z)$ on $|z| = 1 $. Specifically, if $P(z) $ is a polynomial of degree $ n $ with all zeros in $|z| \leq 1 $, then \\
\begin{align}
n\displaystyle\max_{|z|=1} |P(z)| \leq 2 \displaystyle\max_{|z|=1} |P^\prime(z)|.
\end{align}
This inequality is sharp and equality holds for $P(z) = \alpha z^n + \beta $, where $|\alpha| = |\beta| \neq 0 $.\\
\indent Malik \cite{F} extended $(1)$ and proved that if $P(z)$ is a polynomial of degree $ n $ with all zeros in $ |z| \leq k$ where $ k \leq 1$ then \\
\begin{align}
n\max_{|z|=1} |P(z)| \leq (1+k) \max_{|z|=1} |P^\prime(z)|.
\end{align}
Equality in (2) is achieved for $ P(z) = (z + k)^n $ where $ k \leq 1 $.\\
\indent On the other hand, for the class of polynomials $P(z)=a_nz^n+\displaystyle\sum_{j=\mu}^{n}a_{n-j}z^{n-j}$, $ 1 \leq \mu \leq n$ of degree $n$ having all their zeros in $|z| \leq k, k \leq 1,$ Aziz and Shah \cite{D} proved that
\begin{align}
n\displaystyle\max_{|z|=1} |P(z)| \leq (1 + k^{\mu})\displaystyle\max_{|z|=1} |P'(z)|-\frac{n}{k^{n-\mu}}\displaystyle\min_{|z|=k}|P^{\prime}(z)|.
\end{align}
\indent Malik~\cite{G} generalized (1) by replacing its left-hand side by the factor involving the integral mean of \( P(z) \) on \( |z| = 1 \). Infact, he proved that for a polynomial \( P(z) \) of degree \( n \) having all its zeros in \( |z| \leq 1 \),

\begin{align}
n\Big\{\int_{0}^{2\pi}\Big|P(e^{i\theta})|^pd\theta\Big\}^{1/p}\leq\Big\{\int_{0}^{2\pi}|1+e^{i\theta}|^pd\theta\Big\}^{1/p}\displaystyle\max_{|z|=1}|P^{\prime}(z)|, \qquad p>0.
\end{align}
\indent Aziz~\cite{A} generalised $(2)$ by proving that for a polynomial \( P(z) \) of degree \( n \) with all its zeros in \( |z| \leq k \) where \( k \leq 1 \),

\begin{align}
n\Big\{\int_{0}^{2\pi}\Big|P(e^{i\theta})|^pd\theta\Big\}^{1/p}\leq\Big\{\int_{0}^{2\pi}|1+ke^{i\theta}|^pd\theta\Big\}^{1/p}\displaystyle\max_{|z|=1}|P^{\prime}(z)|,\qquad p>0.
\end{align}
Inequality (5) reduces to the inequality (2) by letting $p\rightarrow\infty$. For $k=1$, it reduces to $(4)$.\\
\indent Let $D_{\alpha}P(z)$ denote the polar derivative \cite{H} of a polynomial $P(z)$ of degree $n$ with respect to a complex number $\alpha$, then
$$ D_{\alpha}P(z):= nP(z)+(\alpha-z)P^{\prime} (z). $$
 Note that the polynomial $D_{\alpha}P(z)$ is of degree at most $n-1$ and it generalizes the ordinary derivative $P^{\prime}(z)$ of $P(z)$ in the sense that
\[\underset{\alpha\rightarrow\infty}{\lim}\frac{D_{\alpha}P(z)}{\alpha}=P^{\prime}(z)\] uniformly with respect to $z$ for $\left|  z\right| \leq R,R>0.$\\
\indent Aziz and Rather~\cite{D} extended $(2)$ to the polar derivative and proved that if all the zeros of \( P(z) \) lie in \( |z| \leq k \) with \( k \leq 1 \), then for \( \alpha \in \mathbb{C} \) with \( |\alpha| \geq k \), the result holds

\begin{align}
n(|\alpha|-k)\displaystyle\max_{|z|=1}|P(z)|\leq(1+k)\displaystyle\max_{|z|=1}|D_{\alpha}P(z)|.
\end{align}
For the class of lacunary type polynomials $P(z)=a_nz^n+\displaystyle\sum_{\j=\mu}^{n}a_{n-j}z^{n-\j}, 1 \leq j \leq n,$ of degree $n$ having all their zeros in $|z| \leq k$ where $k \leq 1$, Aziz and Rather \cite{C} also proved that for every complex number $\alpha$ with $|\alpha| \geq k^{\mu}$
\begin{align}
n(|\alpha|-k^{\mu})\displaystyle\max_{|z|=1}|P(z)|\leq(1+k^{\mu})\displaystyle\max_{|z|=1}|D_{\alpha}P(z)|.
\end{align}
Rather and Mir \cite{NM} refined inequality $(7)$ and extended inequality $(3)$ to the polar derivative and proved that if $P(z) = a_n z^n + \displaystyle\sum_{j=\mu}^{n-1} a_{n-j} z^{n-j}$, $1 \leq \mu \leq n$, is a polynomial of degree $n$ having all its zeros in $|z| \leq k$, $k \leq 1$, then for every complex number $\alpha$ with $|\alpha| \geq k^{\mu},$
\begin{align}
\displaystyle n(|\alpha| - k^{\mu})\max_{|z|=1}|P(z)| + \frac{n(|\alpha| - k)}{k^{n - \mu}}\displaystyle\min_{|z|=1}|P(z)| \leq (1 + k^{\mu})\displaystyle\max_{|z|=1}|D_{\alpha}P(z)|.
\end{align}
Recently Roshan Lal et al. \cite{K} proved that if  \( P(z) \) is a polynomial of degree \( n \) with all zeros in \( |z| \leq k \) where \( k \leq 1 \), then for any complex number \( \alpha \) with \( |\alpha| \geq k \) and for each \( p > 0 \), 
\begin{align}
n(|\alpha| - k) \left( \int_0^{2\pi} |P(e^{i\theta})|^p d\theta \right)^{\frac{1}{p}}
\leq
\left( \int_0^{2\pi} |1 + k e^{i\theta}|^p d\theta \right)^{\frac{1}{p}} 
\left( \max_{|z|=1} |D_\alpha P(z)| - \frac{nm}{k^{n-1}} \right)
\end{align}
where $m = \displaystyle\min_{|z|=k} |P(z)|$.\\
In the same paper they also considered class of lacunary type polynomials \( P(z) = a_n z^n + \displaystyle\sum_{j=\mu}^{n} a_{n-j} z^{n-j} \), \( 1 \leq \mu \leq n \) of degree \( n \) with all its zeros in \( |z| \leq k \), \( k \leq 1 \) and proved that for every complex number \( \alpha \) with \( |\alpha| \geq k^\mu \) and for each \( p > 0 \), 
\begin{align}
n(|\alpha| - k^\mu) 
\left( \int_0^{2\pi} |P(e^{i\theta})|^p d\theta \right)^{\frac{1}{p}}
\leq 
\left( \int_0^{2\pi} |1 + k^\mu e^{i\theta}|^p d\theta \right)^{\frac{1}{p}} 
\left( \max_{|z|=1} |D_\alpha P(z)| - \frac{nm}{k^{(n-\mu)}} \right),
\end{align}
where $m = \displaystyle\min_{|z|=k} |P(z)|.$\\
\indent If we let $p\rightarrow\infty$ in (10) and using the fact that $\displaystyle\lim_{p\rightarrow\infty}\Big\{\frac{1}{2\pi}\int_{0}^{2\pi}|P(e^{i\theta}|^pd\theta\Big\}^{1/p}=\displaystyle\max_{|z|=1}|P(z)|$, it does not yield inequality $(8)$. It is natural to seek the integral inequality analogous to inequality (10) which leads to inequality (8) as a special case.
\section{Main Results}
In this paper we present the following integral inequality, which as a special case leads to the inequality (8).
\begin{theorem}\label{th1}
If $P(z) = \displaystyle\sum_{j=0}^{n} a_j z^j$ is a polynomial of degree $n$ having all its zeros in $\left|z\right|\leq k$ where $k\leq 1,$ and $m=\displaystyle \min_{|z|=k}|P(z)|$, then for $\alpha, \beta \in \mathbb C$ with $\left|\alpha\right|\geq k,$  $|\beta|\leq1$   and for each $ p > 0$,
\begin{align*}
n(|\alpha|-k)\Big(\int_{0}^{2\pi}\Big|\frac{P(e^{i\theta})+\frac{m\beta}{k^{n-1}}}{|D_{\alpha}P(e^{i\theta})|-\frac{nm}{k^{n-1}}}\Big|^pd\theta\Big)^{1/p}\leq\Big(\int_{0}^{2\pi}|1+ke^{i\theta}|^pd\theta\Big)^{1/p}.
\end{align*}
\end{theorem}
The result is best possible and equality holds for $P(z)= (z-k)^{n},\quad \alpha>0$.
\begin{remark}\normalfont
Since \( |\alpha| \geq k \), therefore from the inequality $(17)$ with $\mu =1$ (in the proof of Theorem 2), we have for \( |z| = 1 \),
\[
|D_{\alpha}P(z)| - \frac{mn}{k^{n-1}} \geq (|\alpha| - k) |P(z)| \geq 0,
\]
which implies for \( |z| = 1 \),
\[
\big||D_{\alpha}P(z)| - \frac{mn}{k^{n-1}}\big| = |D_{\alpha}P(z)| - \frac{mn}{k^{n-1}} \leq \max_{|z|=1} |D_{\alpha}P(z)| - \frac{mn}{k^{n-1}}.
\]
\end{remark}
Using this fact, the following result follows immediately from Theorem 1. 
\begin{corollary}
If $P(z)$ is a polynomial of degree $n$ having all its zeros in $\left|z\right|\leq k$ where $k\leq 1,$ then for $\alpha, \beta \in \mathbb C$ with $\left|\alpha\right|\geq k,$  $|\beta|\leq1$   and for each $ p > 0$, $m=\displaystyle \min_{|z|=k}|P(z)|$
\begin{align*}
n(|\alpha|-k)\Big(\int_{0}^{2\pi}\Big|P(e^{i\theta})+\frac{m\beta}{k^{n-1}}\Big|^pd\theta\Big)^{1/p}\leq\Big(\int_{0}^{2\pi}|1+ke^{i\theta}|^pd\theta\Big)^{1/p}{\Big(\displaystyle\max_{|z|=1}|D_{\alpha}P(z)|-\frac{nm}{k^{n-1}}}\Big).
\end{align*}
\end{corollary}
Dividing both the sides of Corollary 1 by $|\alpha|$ and letting $|\alpha|\rightarrow \infty$, we get the following result.
\begin{corollary}
If $p(z)$ is a polynomial of degree $n$ having all its zeros in  $|z| \leq k , k \leq 1$, then for $p>0$ and $|\beta| \leq 1$,
\begin{align*}
n\Big(\int_{0}^{2\pi}\Big|P(e^{i\theta}) + \frac{m\beta}{k^{n-1}}\Big|^{p}d\theta\Big)^{1/p}\leq\Big(\int_{0}^{2\pi}|1+ke^{i\theta}|^pd\theta\Big)^{1/p}\displaystyle\max_{|z|=1}|P^{\prime}(z)|.
\end{align*}
\end{corollary}
Instead of proving Theorem 1, we prove the following more general result which includes Theorem 1 as a special case.
\begin{theorem}\label{th2} \normalfont
If $P(z) = a_n z^n + \displaystyle\sum_{j=\mu}^{n} a_{n-j} z^{n-j} $, $1 \leq \mu \leq n$, is a polynomial of degree \(n\) having all its zeros in \( |z| \leq k \), $k \leq 1$ then for $\alpha, \beta \in \mathbb C$ with $\left|\alpha\right|\geq k^\mu,$  $|\beta|\leq1$   and for each $ p > 0$,
\begin{align*}
n(|\alpha|-k^\mu)\Big(\int_{0}^{2\pi}\Big|\frac{P(e^{i\theta})+\frac{m\beta}{k^{n-\mu}}}{|D_{\alpha}P(e^{i\theta})|-\frac{nm}{k^{n-\mu}}}\Big|^pd\theta\Big)^{1/p}\leq\Big(\int_{0}^{2\pi}|1+k^{\mu}e^{i\theta}|^pd\theta\Big)^{1/p}
\end{align*}
where $m=\displaystyle \min_{|z|=k} |P(z)| $.\\
\end{theorem}
As in the case of Corollary 1, we deduce the following result from Theorem 2:
\begin{corollary}\normalfont
If $P(z) = a_n z^n + \displaystyle\sum_{j=\mu}^{n} a_{n-j} z^{n-j}$, $1 \leq \mu \leq n$ is a polynomial of degree \(n\) having all its zeros in \( |z| \leq k \), $k \leq 1$ then for $\alpha, \beta \in \mathbb C$ with $\left|\alpha\right|\geq k^\mu,$  $|\beta|\leq1$ and for each \(p > 0\),
\begin{align}
n(|\alpha|-k^\mu)\Big(\int_{0}^{2\pi}\Big|P(e^{i\theta})+\frac{m\beta}{k^{n-\mu}}\Big|^pd\theta\Big)^{1/p}\leq\Big(\int_{0}^{2\pi}|1+k^{\mu}e^{i\theta}|^pd\theta\Big)^{1/p}{\Big(\displaystyle\max_{|z|=1}|D_{\alpha}P(z)|-\frac{nm}{k^{n-\mu}}}\Big),
\end{align}
which is the desired integral analogue of inequality (10). If we let $p\rightarrow \infty$ in inequality $(11)$, we get the inequality $(8)$.
\end{corollary}
For $\beta =0$, we get inequality (10).
Next by using Holder's Inequality, we establish the following result.
\begin{theorem}\label{th3} \normalfont
If $P(z) = a_n z^n + \displaystyle\sum_{j=\mu}^{n-j} a_j z^j $, $1 \leq \mu \leq n$ is a polynomial of degree \(n\) having all its zeros in \( |z| \leq k \) where  $k \leq 1$, then for complex numbers \(\alpha\) and $\beta$ with \(|\alpha| \geq k^\mu\), $|\beta| \leq 1$ and for each \(r > 1\), \(s > 1\) with $\frac{1}{r}+\frac{1}{s}=1$,
\begin{align}
n(|\alpha|-k^\mu)\Big\{\int_{0}^{2\pi}\Big|P(e^{i\theta})+\frac{m\beta}{k^{n-\mu}}\Big|^pd\theta\Big\}^{1/p}\leq\Big\{\int_{0}^{2\pi}|1+k^{\mu}e^{i\theta}|^{pr}d\theta\Big\}^{1/pr}{\Big\{\int_{0}^{2\pi}\Big(|D_{\alpha}P(e^i{\theta})|-\frac{nm}{k^{n-\mu}}}\Big)^{ps}d\theta\Big\}^{1/ps}.
\end{align}
where $m=\displaystyle \min_{|z|=k}|P(z)|$.
\end{theorem}
\begin{remark} \normalfont
Letting $s \rightarrow \infty$ in (12) and noting that $r \rightarrow 1$, we get Theorem 2.
\end{remark}
\section{Lemma}
For the proof of above theorems, we need the following lemma, which is due to N. A. Rather \cite{J}.
\begin{lemma}\label{le3}
Let $P(z)=a_nz^n+\displaystyle\sum_{j=\mu}^{n}a_{n-j}z^{n-j}$, $ 1 \leq \mu \leq n $, be a polynomial of degree n having all its zeros in $\left|z\right|\leq k$ where $k\leq 1,$ then for $|z|=1$, 
\begin{align*}
k^{\mu}|P'(z)| \geq |Q'(z)|+\frac{nm}{k^{n-\mu}}
\end{align*}
\end{lemma}
where $Q(z)=z^n\overline{P(\frac{1}{\overline{z}})}$ and $m = min_{|z|=1|}|P(z)|.$
\section{Proof of Theorems}
In this section we prove Theorem 2 and the proof of Theorem 1 follows from Theorem 2 by taking $\mu=1$.
\begin{proof}[\bf{ Proof of Theorem \ref{th2}}] \normalfont
 Since $Q(z)=z^n\overline{P(\frac{1}{\overline{z}})}$, therefore $P(z)=z^n\overline{Q(\frac{1}{\overline{z}})}$ and it can be easily verified that 
 \begin{align}
 |Q'(z)| = |nP(z) - zP'(z)| \qquad \text{and} \qquad |P'(z)|=|nQ(z)-zQ'(z)|,   \qquad for \quad |z|=1.
 \end{align}
By Lemma 1, we have for every $\beta$ with $|\beta| \leq 1$ and $|z|=1$,
\begin{align*}
\Big|Q'(z)+\frac{nm\overline{\beta}z^{n-1}}{k^{n-\mu}}\Big| &\leq |Q'(z)|+\frac{nm|\beta|}{k^{n-\mu}}\\
                                                   &\leq |Q'(z)|+\frac{nm}{k^{n-\mu}}\\
                                                  &\leq k^{\mu}|P'(z)|\\
                                                  &=k^\mu|nQ(z)-zQ'(z)| \qquad |z|=1.
\end{align*}
That is,
\begin{align}
\Big|Q'(z)+\frac{nm\overline{\beta}z^{n-1}}{k^{n-\mu}}\Big|\leq k^\mu|nQ(z)-zQ'(z)| \qquad |z|=1.
\end{align}
Since $P(z)$ has all its zeros in $|z| \leq k$ where $k \leq 1$, by Gauss-Lucas Theorem that all the zeros of $P^\prime(z)$ also lie in $|z| \leq k$, $k\leq1$. Hence the polynomial 
\begin{align*}
z^{n-1}\overline{P^\prime(\frac{1}{\overline{z}})}=nQ(z)-zQ^{\prime}(z)
\end{align*}
has all its zeros in $|z|\geq1$. Therefore, it follows from $(3)$ that the function
\begin{align*}
w(z)=\frac{z(Q^{\prime}(z)+\frac{nm\overline{\beta}z^{n-1}}{k^{n-\mu}})}{k^\mu(nQ(z)-zQ^{\prime}(z))}
\end{align*}
is analytic in $|z|\leq1$ and $|w(z)|\leq1$ for $|z|=1$. Also $w(0)=0$. Thus the function $1+k^{\mu}w(z)$ is subordinate to the function $1+k^{\mu}z$ for $|z|\leq1$. Hence by well known property of subordination [\cite{E}, P.422] for each $p>0$,
\begin{align}
\int_{0}^{2\pi}\Big|1 + k^{\nu}w(e^{i\theta})\Big|^pd\theta \leq \int_{0}^{2\pi}|1+k^{\mu}e^{i\theta}|^{pr}d\theta.
\end{align}
Now
\begin{align*}
|1+k^{\mu}w(z)|=\frac{n\Big|Q(z)+\frac{m\overline{\beta}z^n}{k^{n-\mu}}\Big|}{|nQ(z)-zQ^{\prime}(z)|}.\\
\end{align*}
This gives with the help of $(13)$,
\begin{align}
|1+k^{\mu}w(z)|=\frac{n\Big|Q(z)+\frac{m\overline{\beta}z^n}{k^{n-\mu}}\Big|}{|P^{\prime}(z)|} \qquad for \qquad |z|=1.
\end{align}
Also by $(13)$ for $|z|=1$, we have 
\begin{align*}
 |D_{\alpha}P(z)|&=|nP(z)+(\alpha-z)P^{\prime} (z)|\\
 &= |nP(z) - zP^{\prime}(z) + \alpha P^{\prime}(z)|\\
 &\geq |\alpha||P^{\prime}| - |nP(z0- zP^{\prime}(z)|\\
 &= |\alpha||P^{\prime}| - |Q^{\prime}(z)|.
 \end{align*}
 Using Lemma \ref{le3}, we get for $|z|=1$,
 \begin{align*}
 |D_{\alpha}P(z)| &\geq|\alpha||P^{\prime}(z)|-k^\mu|P^{\prime}(z)|+\frac{nm}{k^{n-\mu}}
 \end{align*}
 \begin{align}
 &=(|\alpha|-k^\mu)|P^{\prime}(z)|+\frac{nm}{k^{n-\mu}}.
\end{align}
Which gives
\begin{align}
|P^{\prime}(z)|\leq\frac{ |D_{\alpha}P(z)|-\frac{nm}{k^{n-\mu}}}{(|\alpha|-k^\mu)},\quad |z|=1
\end{align}
Together $(16)$ and $(18)$ leads to
\begin{align}
|1+k^{\mu}w(z)| \geq (|\alpha|-k^\mu)\frac{n\Big|Q(z)+\frac{m\overline{\beta}z^n}{k^{n-\mu}}\Big|}{ |D_{\alpha}P(z)|-\frac{nm}{k^{n-\mu}}}, \quad \text{for} \quad |z|=1.
\end{align}
Using $(19)$ in $(15)$ we obtain for each $p>0$,
\begin{align*}
n(|\alpha|-k^\mu)\Big(\int_{0}^{2\pi}\Big|\frac{Q(e^{i\theta})+\frac{m\bar{\beta} e^{in\theta}}{k^{n-\mu}}}{|D_{\alpha}P(e^{i\theta})|-\frac{nm}{k^{n-\mu}}}\Big|^pd\theta\Big)^{1/p}\leq\Big(\int_{0}^{2\pi}|1+k^{\mu}e^{i\theta}|^pd\theta\Big)^{1/p},
\end{align*}
that is,
\begin{align*}
n^{p}(|\alpha|-k^\mu)^{p}\int_{0}^{2\pi}\Big|\frac{Q(e^{i\theta})+\frac{m\bar{\beta} e^{in\theta}}{k^{n-\mu}}}{|D_{\alpha}P(e^{i\theta})|-\frac{nm}{k^{n-\mu}}}\Big|^pd\theta \leq\int_{0}^{2\pi}|1+k^{\mu}e^{i\theta}|^pd\theta, 
\end{align*}
Equivalently for each $p>0$, we have 
\begin{align*}
n(|\alpha|-k^\mu)\Big(\int_{0}^{2\pi}\Big|\frac{P(e^{i\theta})+\frac{m\beta}{k^{n-\mu}}}{|D_{\alpha}P(e^{i\theta})|-\frac{mn}{k^{n-\mu}}}\Big|^pd\theta\Big)^{\frac{1}{p}} \leq \Big(\int_{0}^{2\pi}|1+k^{\mu}e^{i\theta}|^pd\theta\}\Big)^{\frac{1}{p}}.
\end{align*}
This completes the proof of Theorem $2$.
\end{proof}
\begin{proof}[\bf{ Proof of Theorem \ref{th3}}] \normalfont
 Since all the zeros of polynomial $P(z)=a_nz^n+\displaystyle\sum_{j=\mu}^{n}a_{n-j}z^{n-j}$, $ 1 \leq \mu \leq n $, lie in $|z| \leq k$, $k\leq 1$, therefore, proceeding as in the proof of Theorem 2, we have from inequality $(19)$ for each \( p > 0 \),
\begin{equation}
n^p (|\alpha| - k^{\mu})^p \int_0^{2\pi}  \big|nQ(e^{i\theta}) + \frac{m\overline{\beta}e^{in\theta}}{k^{n-\mu}}\big|^{p} \, d\theta 
\leq \int_0^{2\pi} |1 + k^{\mu} w(e^{i\theta})|^p \left( |D_{\alpha}P(e^{i\theta})| - \frac{mn}{k^{n-\mu}} \right)^{p} \, d\theta.
\end{equation}
Using Hölder's inequality and $(15)$ for \( r > 1, s > 1 \) with \( \frac{1}{r} + \frac{1}{s} = 1 \), we get

\begin{align*}
\begin{aligned}
n^p (|\alpha| - k^{\mu})^p \int_0^{2\pi} \left|ne^{in\theta}\overline{P(e^{i\theta})} + \frac{m\overline{\beta}e^{in\theta}}{k^{n-\mu}}\right|^{p} d\theta 
&\leq \left\{ \int_0^{2\pi} |1 + k^{\mu} w(e^{i\theta})|^{pr} d\theta \right\}^{1/r} 
\left\{ \int_0^{2\pi} \left(|D_{\alpha}P(e^{i\theta})| - \frac{mn}{k^{n-\mu}} \right)^{ps} d\theta \right\}^{1/s}
\end{aligned}
\end{align*}
\begin{align*}
\qquad \leq \left\{ \int_0^{2\pi} |1 + k^{\mu}e^{i\theta}|^{pr} d\theta \right\}^{1/r} \left\{ \int_0^{2\pi} \left(|D_{\alpha}P(e^{i\theta})| - \frac{mn}{k^{n-\mu}} \right)^{ps} d\theta \right\}^{1/s}.
\end{align*}
That is,
\begin{align*}
\begin{aligned}
n (|\alpha| - k^{\mu})\left\lbrace \int_0^{2\pi} \left|nP(e^{i\theta}) + \frac{m\beta}{k^{n-\mu}}\right|^{p} d\theta\right\rbrace^{\frac{1}{p}}
&\leq \left\{ \int_0^{2\pi} |1 + k^{\mu}e^{i\theta}|^{pr} d\theta \right\}^{1/pr} \left\{ \int_0^{2\pi} \left(|D_{\alpha}P(e^{i\theta})| - \frac{mn}{k^{n-\mu}} \right)^{ps} d\theta \right\}^{1/ps}.
\end{aligned}
\end{align*}
This completes the proof of Theorem 3.
\end{proof}
\begin{proof}[\bf{ Proof of Theorem \ref{th1}}] \normalfont
The proof of Theorem $1$ follows along the same lines by putting $\mu=1$, in Theorem $2$.
\end{proof}


\begin{thebibliography}{99}

\bibitem{A} A. Aziz, Integral mean estimates for polynomials with resticted zeros, J. Approx. Theory, 55 (1988),
232-239
\bibitem{B} A. Aziz and N. A. Rather, A refinement of a theorem of Paul Turan concerning
polynomials. Math Ineq. Appl., 1(1998), 231-238.
\bibitem{C} A. Aziz and N. A. Rather, Inequalities for the polar derivative of a polynomial
with restricted zeros, Math. Balk., 17(2003), 15-28.
\bibitem{D}  A. Aziz and W. M. Shah, Integral mean estimates for polynomials, Ind. J. Pure Appl. Math., 28(1997), 1413-1420.
\bibitem{GV} R. B. Gardner, N. K. Govil and G. V. Milovanović, Extremal problems and inequalities of Markov-Bernstein type for algebraic polynomials. Academic Press, 2022.
\bibitem{E} E. Hille, Analytic Function Theory II Ed, Ginn and Company, New York, (1962). 
\bibitem{F}  M. A. Malik, On the derivative of a polynomial, J. Lond. Math. Soc., Second Series 1(1969), 57-60
\bibitem{G} M. A. Malik, An integral mean estimates for polynomials, Proc. Amer. Math. Soc.,
91(1984), 281-284.
\bibitem{H} M. Marden, Geometry of Polynomials. Second Edition. Mathematical Surveys,
No. 3 American Mathematical Society, Providence, R.I. 1966 
\bibitem{I} G. V. Milvonic, D. S. Mitrovonic and Th. M. Rassias, Topics in Polynomials: Extremal Problem, Inequalities, Zeros, World ScientiÖc, Singapore, 1994.
\bibitem{J} N. A. Rather, Extremal properties and Location of the zeros of polynomials, Ph.D.
Thesis, University of Kashmir, 1998.
\bibitem{NM} N. A. Rather and M. I. Mir, Some refinements of inequalities for the polar derivative of polynomials with restricted zeros, Int. J. Pure and Appl. Math., 41(2007),
1065-1074
\bibitem{K} Roshan Lal Keshtwal, Susheel Kumar and Imran Ali, Generalized $L_p$ Inequalities For The Polar Derivative Of Polynomial. J\~n\=an\=abha, 54 (1) (2024), 83-87.
\bibitem{TU} P. Turán, Über die Ableitung von Polynomen,'' \textit{Compos. Math.}, 7(1939),  89--95.


\end{thebibliography}
\end{document}